\documentclass[11pt]{amsart}
\usepackage{amsmath}
\usepackage{amsthm}
\usepackage{amssymb}
\usepackage{romannum}
\newtheorem{thm}{Theorem}[section]
\newtheorem{lem}[thm]{Lemma}

\newtheorem{cor}[thm]{Corollary}

\newtheorem{rem}[thm]{Remark}

\usepackage[margin=1in]{geometry}

\begin{document}

\title{Additive twist of $- (\zeta '/ \zeta)(s)$, nontrivial zeros of $L(s, \chi)$,
and Continuity}
\author{Hisanobu Shinya}
\address{Toyama prefecture, Japan. The author has no academic affiliation.}
\email{zhilijohntaro683@gmail.com}

\begin{abstract}
Suppose that the Riemann hypothesis is false and let $\rho_{*} = 1/2 + \eta_{*} + i \gamma_{*}$
($\eta_{*} > 0$) be a nontrivial zero of the Riemann $\zeta$-function off the critical line.
We establish an asymptotic
relation (as $\gamma_{*} \to \infty$) linking the residues of the exponential Dirichlet series
$
\sum_{n \geq 1} \Lambda(n) e^{- 2\pi i p n } n^{-s}$
at the corresponding nontrivial zeros of associated Dirichlet $L$-functions 
to an explicit function that is continuous in $p$.
This asymptotic relation holds for all rational numbers $p = a/b \in (0, 1)$ with
$b \ll \log \gamma_{*}$.
We discuss the implications of this continuity for local distribution of nontrivial zeros of Dirichlet $L$-functions.
\end{abstract}

\maketitle

\section{Introduction}

Ever since the Riemann hypothesis was proposed in 1859, researchers throughout centuries have never
quit trying to resolve the problem, and at times they have succeeded in inventing new useful machineries to
analyze it. In this paper, let us invent another one. 

Several results of statistical type such as $N(\sigma, T) \ll T^{3(1 - \sigma)/(2- \sigma)} (\log T)^{5}$
\cite{T}, where
$N(\sigma, T)$ denotes the number of the zeros of the Riemann $\zeta$-function in the region
$\{ a + bi:  a  > \sigma, 0 < b \leq T \}$,
have been established. The main goal of this paper is to give a local-type result that suggests a new phenomenon
characteristic of local zeros of Dirichlet $L$-functions $L(s, \chi)$ lying near a nontrivial zero
$\rho_{n_{1}, p}$ of the
Riemann $\zeta$-function with $\text{Re}(\rho_{n_{1}, p}) > 1/2$. 

In our method, the main tool is the function
\[
M(s, p) \equiv
\sum_{n \geq 1} \Lambda(n) e^{- 2\pi i p n } n^{-s},
\]
where $\Lambda(n)$ is the arithmetical Mangoldt $\Lambda$-function and
$p \in \mathbb{Q} \cap (0, 1)$.
With modern terminology,
the function $M(s, p)$ is recognized as the additive twist
of the Dirichlet series
\[
-\frac{\zeta'}{\zeta}(s) = \sum_{n \geq 1} \frac{\Lambda(n)}{n^{s}}.
\]

The use of additive twists is found in the literature; for instance, Jutila
made use of the twisted series
\[
\sum_{n \geq 1} \frac{d(n)e(np)}{n^{s}},
\]
where $d(n)$ denotes the number of divisors of $n$,
to analyze the fourth moment of the Riemann $\zeta$-function
on the critical line \cite{008}.
Also, the partial sum
\begin{equation} \label{finiteone}
S(\alpha) \equiv \sum_{k \leq N} \Lambda(k) e^{-2 \pi i k \alpha}
\end{equation}
plays the key role in the study of weak Goldbach conjecture \cite{H}.
One appearance of the series $M(s, p)$ in the literature
is found in \cite{HL}, where Hardy and Littlewood considered, instead of (\ref{finiteone}),
the infinite sum
\[
\sum_{(q, n) = 1} \Lambda(n)x^{n}
\]
in their study of the Goldbach problem. 

In our study, we too shall use the function $M(s, p)$.
It becomes especially useful when rewritten as the following lemma.
The identity appears in \cite{HL}, but we give its proof here.
\begin{lem} \label{rwri} 
Let $0 < a/b < 1$ be a rational number. For all $\text{Re}(s) > 0$, we have
\[
M(s, a/b) -  \sum_{q|b} \sum_{l \geq 1} (\log q) q^{-sl} e^{- 2 \pi i a q^{l}/b}
 = -  \sum_{\chi}  A(a, b; \chi) \frac{L'}{L}(s, \chi),
\]
where the sum $\sum_{q|b}$ runs through all the prime divisors of $b$,
$L(s, \chi)$ is the Dirichlet $L$-function with a character $\chi$ modulo $b$, and 
\[
A(a, b; \chi) \equiv
\frac{\sum_{1 \leq x \leq b} \overline{\chi}(x) e^{- 2\pi i a x / b}   }{\phi(b)}.
\]
In addition, we have for $\text{Re}(s) > 0$,
\begin{equation} \label{ganbaruman}
\begin{split}
E_{M}(s, a/b) 
&\equiv  \sum_{q|b} \sum_{l \geq 1} (\log q) q^{-sl} e^{- 2 \pi i a q^{l}/b} \\
& = \sum_{q|b} (\log q)
\left[\sum_{l = 1}^{k - 1}q^{-sl} e^{- 2 \pi i a q^{l}/b}
+ \frac{1}{1 - q^{-ds}} \sum_{r = 0}^{d - 1} q^{-s(k + r)}e^{- 2 \pi i a q^{k + r}/b}
\right],
\end{split}
\end{equation}
where $k$, $q$, and $m$ satisfy $(q, m) = 1$ and
\[
q^{k}m = b,
\]
and $d$ is the least positive integer satisfying
\[
q^{d} \equiv 1 \pmod{m}.
\]
Hence, the function
$E_{M}(s, a/b)$
can be extended meromorphically to the whole complex plane. 
\end{lem}

\begin{proof}
The first part of the lemma follows easily by applying the formula, for $(k, b) = 1$,
\begin{equation} \label{orthprope}
\begin{split}
e^{- 2\pi i a k /b} 
&
=
\sum_{1 \leq x \leq b}e^{-2\pi i a x / b}
\sum_{\chi:(\mathbb{Z} / b)^{*} \to \mathbb{C}^{*}} \frac{\chi(k)\overline{\chi}(x)}{\phi(b)}
\\
& = \sum_{\chi} \chi(k)
\frac{\sum_{1 \leq x \leq b} \overline{\chi}(x) e^{- 2 \pi i a x / b}   }{\phi(b)} \\
& =  \sum_{\chi} \chi(k) A(a, b; \chi),
\end{split}
\end{equation}
which can be easily shown with the well-known identity concerning
Dirichlet characters valid for $(a, m) = 1$ \cite{H}
\[
\frac{\sum_{\chi \in (\mathbb{Z} / m)^{*}} \overline{\chi}(x) \chi(a)  }{\phi(m)}
= 
\begin{cases}
1   &  x \equiv a \text{(mod $m$)}, \\
0   &  \text{otherwise};
\end{cases}
\]
here, for each $k$,
there exists only one $1 \leq k' \leq b$
which is equal to $k$ modulo $b$
so that the sum of (\ref{orthprope}) is equal to $e^{- 2\pi i a k' / b} = e^{- 2\pi i a k / b}$.

We rewrite
\begin{equation} \label{MtoL}
\begin{split}
-  \sum_{\chi}  A(a, b; \chi) \frac{L'}{L}(s, \chi)
& =  \sum_{\chi}  A(a, b; \chi)  \sum_{n \geq 1} 
\frac{
\Lambda(n)  \chi(n)
}{n^{s}} \\
&  = \sum_{n \geq 1} \frac{\Lambda(n)}{n^{s}}
 \sum_{1 \leq x \leq b} e^{-2 \pi i a x/b} \sum_{\chi} \frac{\overline{\chi}(x) \chi(n)}{\phi(b)}.
\end{split}
\end{equation}

By the fact that the formula (\ref{orthprope}) is applicable to all $n$'s which are powers of primes other than 
powers of primes $q$ satisfying $q|b$,
we have
\[
M(s, a/b) - \sum_{q|b} \sum_{l \geq 1} (\log q) q^{-sl} e^{- 2 \pi i a q^{l}/b}
 = -  \sum_{\chi}  A(a, b; \chi) \frac{L'}{L}(s, \chi).
\]

This completes the proof of the first part of the lemma.

We give a proof of the second part as follows.\footnote{The argument is due to the 
large language model AI Gemini.} 
By the definitions of $d$, we have
\[
q^{l + d} - q^{l} = q^{l}(q^{d} - 1) \equiv 0 \pmod{q^{k} m = b}
\]
whenever $k \leq l$, or
\[
q^{l + d} \equiv q^{l} \pmod{b}.
\]
This implies that the function $e^{- 2 \pi i a q^{l}/b}$ is periodic in $l$
with period $d$.

We rewrite with $l = k + r + jd$,  $0 \leq r \leq d - 1$, $j \geq 0$,
\[
\begin{split}
\sum_{l \geq 1} q^{-sl} e^{- 2 \pi i a q^{l}/b}
& = \sum_{l = 1}^{k - 1}  q^{-sl} e^{- 2 \pi i a q^{l}/b}
+ \sum_{l \geq k}  q^{-sl} e^{- 2 \pi i a q^{l}/b} \\
& =  \sum_{l = 1}^{k - 1}  q^{-sl} e^{- 2 \pi i a q^{l}/b}
+ \sum_{r = 0}^{d - 1} \sum_{j \geq 0}  q^{-s(k + r + jd)} e^{- 2 \pi i a q^{k + r + jd}/b}.
\end{split}
\]
By the periodicity of  $e^{- 2 \pi i a q^{l}/b}$, the second sum becomes
\[
\begin{split}
 \sum_{r = 0}^{d - 1} \sum_{j \geq 0}  q^{-s(k + r + jd)} e^{- 2 \pi i a q^{k + r + jd}/b} 
&  =  \sum_{r = 0}^{d - 1} \sum_{j \geq 0}  q^{-s(k + r + jd)} e^{- 2 \pi i a q^{k + r}/b} \\
& = \frac{1}{1 - q^{-ds}}\sum_{r = 0}^{d - 1}  q^{-s(k + r)} e^{- 2 \pi i a q^{k + r}/b}.
\end{split}
\]

This completes the proof of the lemma.

\end{proof}

To achieve our goal, we base our argument on 
the following theorem.

\begin{thm} \label{main1}
Let $0 < c <1$, $\delta > 0$, $c + \delta + \kappa > 1$, $c' > 0$, $\kappa - c' > 1$, $\delta + c' > 1$, 
and $p$ be any real numbers satisfying $0 < p < 1$.
Then we have
\begin{equation} \label{thm1}
\begin{split}
& \frac{1}{2\pi i}\int_{(c)} 
\sum_{n \geq 2} \frac{\Lambda(n) e^{-2 \pi i p n} }{
n^{s + \delta + \kappa}}
\Gamma(s) \Gamma(s+ \delta)\Gamma(1-s) (2 p\pi)^{-s}e^{-\pi i s / 2}ds \\
& = \Gamma(1 + \delta)
(- \frac{\zeta'}{\zeta}(\kappa))
[ (2 \pi p)^{-\delta}  e^{\delta \pi i / 2}\Gamma( - \delta)
+ K(p, 1 + \delta)]
\\
& -  \frac{\Gamma(1 + \delta)}{2 \pi i }
\int_{c' - i \infty}^{c' + i\infty}
- \frac{\zeta'}{\zeta}
(\kappa - w) \\
& \times 
 [(2 \pi p)^{\delta + w}  e^{(\delta + w ) \pi i / 2}\Gamma( - \delta - w)
+ K(p, 1 + \delta + w)]
\frac{dw}{w},
\end{split}
\end{equation}
where $K(p, \xi)$ denotes a function defined for $0 < \text{Re}(\xi) < 1$ by
\[
K(p, \xi) \equiv - \int_{0}^{1} e^{- 2 \pi i p r } r^{-\xi} dr.
\]
\end{thm}

The proof of Theorem \ref{main1} will be given in the next section.

From Theorem \ref{main1}, it is easy to deduce the following theorem.
\begin{thm} \label{supraclaim}
Let $h > 0$, $\delta = \nu + iT$ ($\nu > 0$) and $p = a/b < 1$. We have
\begin{equation} \label{supraidentity}
\begin{split}
& -\sum_{\chi }\sum_{\rho_{\chi}} A(a, b; \chi) (2 p \pi)^{- \rho_{\chi} + \delta}
e^{-\pi i (\rho_{\chi} - \delta ) / 2}  \Gamma(
\rho_{\chi} - \delta
)
\Gamma(\rho_{\chi})\Gamma(1- \rho_{\chi} + \delta) \\
& + M(0, a/b) \Gamma(\delta) 
+ E_{1} + E_{2} + E_{3} + E_{4}  \\
& =
\sum_{\rho_{j}}
\rho_{j}^{-1}
J(p, 1 + \delta - \rho_{j}) 
+ E_{5} + E_{6},
\end{split}
\end{equation} 
where $\rho_{j}$ denotes the nontrivial zeros of the Riemann $\zeta$-function, $\rho_{\chi}$ those of $L(s, \chi)$,
and
\begin{align}  \label{errors}
E_{1} & \equiv  - \sum_{\chi}A(a, b; \chi) \text{res}_{s = - \delta}\left[ \frac{L'}{L}(s + \delta) 
\Gamma(s) \Gamma(s+ \delta) \Gamma(1-s) (2 p\pi)^{-s}e^{-\pi i s / 2}\right] \notag \\
E_{2} & \equiv - \sum_{\chi^{*}} \sum_{n}A(a, b; \chi^{*}) 
 \Gamma( h^{*}(n) - \delta) \Gamma(h^{*}(n)) 
\Gamma(1-h^{*}(n) + \delta) (2 p\pi)^{- h^{*}(n) + \delta}e^{-\pi i ( h^{*}(n) - \delta) / 2} \notag \\
E_{3} & \equiv  \text{res}_{s = - \delta}[E_{M}(s + \delta, a/b)\Gamma(s) \Gamma(s+ \delta)
 \Gamma(1-s) (2 p\pi)^{-s}e^{-\pi i s / 2} ] \notag \\
E_{4} & \equiv \frac{1}{2\pi i}\int_{(-1/2)} 
M(s + \delta, a/b)
\Gamma(s) \Gamma(s+ \delta)\Gamma(1-s) (2 p\pi)^{-s}e^{-\pi i s / 2}ds \notag \\
E_{5} & \equiv \sum_{\rho_{j}}
\rho_{j}^{-1} K(p,  1 + \delta - \rho_{j}) \notag \\
E_{6} & \equiv - 
[
J(p, \delta) + K(p, \delta)
]   
+  \frac{1}{2 \pi i } \int_{2 - h - i \infty}^{2 - h + i\infty}
- \frac{\zeta'}{\zeta}
( - w) 
[J(p, 1 + \delta + w) + K(p, 1 + \delta + w)]\frac{dw}{w}. \notag 
\end{align}
Here, $\chi^{*}$ denotes imprimitive characters modulo $b$, and
$h^{*}(n)$ the roots on the vertical line $\text{Re}(s) = 0$
of 
\[
L(s, \chi^{*}) = L(s, \chi_{1}) \prod_{p | b} [1 - \chi_{1}(p)p^{-s}],
\]
where $\chi_{1}$ is the primitive character inducing $\chi^{*}$ \cite{H}.

In addition, all the terms on the right side are continuous in $p$.
\end{thm}

The proof of Theorem \ref{supraclaim} will be given in Section 3.

The following is a corollary to Theorem \ref{supraclaim},
which states a consequence of Theorem \ref{supraclaim} for small modulus $b$.

\begin{cor} \label{claim}
Suppose that $M(s, p)$ has a pole at  $s = \rho_{n_{1}, p} = 1/2 + \eta_{n_{1}, p} + i \gamma_{n_{1}, p}$
with $\eta_{n_{1}, p} > 0$ and $p = a/b < 1$ is any rational number satisfying $b \ll \log \gamma_{n_{1}, p}$.
Let $\nu > 0$ be arbitrary and
\[
R_{n_{1}, p} \equiv \{a + bi \in \mathbb{C}: 1/2 \leq a \leq 1, 
|b| \leq A \log \gamma_{n_{1}, p}  \}
\]
for some absolute constant $A$ which is independent of $\rho_{n_{1}, p}$.
Suppose also that
$\eta_{n_{1}, p}$ is the largest number among other such numbers
$\text{Re}(\rho_{\chi}) - 1/2$
associated with other
possible zeros of $L(s, \chi)$'s off the criticial line, at which
$M(s + \nu + i\gamma_{n_{1},p}, p)$ may have poles in $R_{n_{1}, p}$.
Then, we have as $\gamma_{n_{1}, p} \to \infty$,
\begin{equation} \label{mode11}
\begin{split}
& 
\sum_{\substack{
\text{Re}(\rho_{\chi}) = 1/2 + \eta_{n_{1}, p} \\
 \rho_{\chi} - \nu - i\gamma_{n_{1}, p}
 \in R_{n_{1}, p} }}A(a, b; \chi) (2 p \pi)^{- \rho_{\chi} + \delta}
e^{-\pi i (\rho_{\chi} - \delta ) / 2}  \Gamma(
\rho_{\chi} - \delta
)
 e^{\pi \gamma_{n_{1}, p} / 2}\Gamma(\rho_{\chi})\Gamma(1- \rho_{\chi} + \delta) \\
& \sim
\Gamma(1 + \nu + i \gamma_{n_{1}, p})e^{\pi \gamma_{n_{1}, p}/ 2}   \\
& \times 
\sum_{\rho_{j}}
\rho_{j}^{-1}
 (2 \pi p)^{\nu + i \gamma_{n_{1}, p} - \rho_{j} }  e^{(\nu + i \gamma_{n_{1}, p} - \rho_{j}  ) \pi i / 2}\Gamma(
- \nu - i \gamma_{n_{1}, p} + \rho_{j}).
\end{split}
\end{equation}
\end{cor}

\begin{proof}
We multiply both sides of (\ref{supraidentity})
by $e^{\pi \gamma_{n_{1}, p} / 2}$
so that the main terms have the order of magnitude $\gamma_{n_{1}, p}^{c}$.

We note that \cite{H} through $m \equiv ah \pmod b$,
\[
\begin{split}
\frac{1}{\phi(b)}\chi(a) \sum_{m = 1}^{b}\overline{\chi}(m)e_{b}(-m)
& = \frac{1}{\phi(b)} \sum_{m = 1}^{b} \chi(a) \overline{\chi}(m)e_{b}(-m) \\
& = \frac{1}{\phi(b)} \sum_{h = 1}^{b}  \overline{\chi}(h)e_{b}(-ah) \\
& = A(a, b; \chi),
\end{split}
\]
or with the Gaussian sum $\tau(\chi)$,
\[
|A(a, b; \chi)| = \left| \frac{\chi(a) \overline{\tau(\chi)}}{\phi(b)} \right|
= \frac{b^{1/2}}{\phi(b)}.
\]

Thus, by the restriction $b \ll \log \gamma_{n_{1}, p}$,
the terms of the left side of (\ref{mode11}) associated with the roots of Dirichlet $L$-functions
keep being major ones.

Besides, by the condition $b \ll \log \gamma_{n_{1}, p}$, it is easy to 
show that all the terms $E's$ in (\ref{supraidentity}) are minor compared to the left members of (\ref{mode11}).

The leftover $\sum_{\not \in R_{n_{1}, p}}$ is minor, because
its terms decay by the order
\[
\begin{split}
& e^{-\pi i (\rho_{\chi} - \delta ) / 2}  \Gamma(
\rho_{\chi} - \delta
)
\Gamma(1- \rho_{\chi} + \delta)
\Gamma(\rho_{\chi}) e^{\pi \gamma_{n_{1}, p} / 2} \\
& \ll e^{\pi(\gamma_{\chi} - \gamma_{n_{1}, p})/2} e^{- \pi |\gamma_{\chi} - \gamma_{n_{1}, p}|}
|\gamma_{\chi}|^{\eta_{\chi}} e^{- \pi |\gamma_{\chi}| / 2} e^{\pi \gamma_{n_{1}, p} / 2} \\
&\ll  e^{- \pi |\gamma_{\chi} - \gamma_{n_{1}, p}|} |\gamma_{\chi}|^{\eta_{\chi}}.
\end{split}
\]

This completes the proof of the corollary.
\end{proof}

\begin{rem}
A simple example of the corollary is when $p = 1/2, 2/3,$ etc.
``Up to what values of $b$, the residue terms would keep being major?'' is a qeustion of curiosity.
Corollary \ref{claim} is in effect whenever the negation of the Riemann hypothesis is assumed. 
This is because $M(s, a/b)$ expressed as in Lemma \ref{rwri}
always contains a term associated with the principal character
$\chi_{0}$.
Note that $L(s, \chi_{0}) = \zeta(s) \prod_{p|b}(1 - p^{-s})$ for modulo $b$. 
\end{rem}

In Section 4, we will discuss a consequence of Theorem \ref{supraclaim}
for large modulus $b$. As we consider the theorem with $p = a/b$ in
a short interval, the modulus $b$ needs to become large.
The continuity property mentioned in the theorem
gives us a foundational principle governing the local distribution of nontrivial zeros of $L(s, \chi)$'s.

The rest of the article is devoted to the proofs of theorems.

Our argument depends fully on
the following estimate for the $\Gamma$-function \cite{L}
\begin{equation} \label{gf}
\begin{split}
\Gamma(\sigma + it) 
\asymp |\sigma + it|^{\sigma - 1/2} e^{- \pi |t| / 2}
\end{split}
\end{equation}
as $|t| \to \infty$.

\section{Proofs of Theorem \ref{main1}}

In this section, we give two proofs of Theorem \ref{main1} to make the theorem solid.

\subsection{The first proof}

We define \cite{L}
\begin{equation} \label{reviewer1}
I_{1}(x) \equiv \frac{1}{2\pi i}\int_{(c)}  \Gamma(s)  x^{-s} ds = e^{-x}
\end{equation}
and  
\[
I_{2}(x) \equiv \frac{1}{2\pi i}\int_{(c)} \Gamma(s+q-1)\Gamma(1-s) x^{-s} ds,
\]
where $0 < q < 1$ with $1 - q < c < 1$ (these choices are for putting $I_{2}(x)$ in a residue-free form).
  
Shifting the path in $I_{2}$ to the left, we find out that
\begin{equation} \label{reviewer2}
I_{2}(x) = 
\sum_{n=0}^{\infty} \frac{(-1)^{n} \Gamma(n+ q) x^{n+q-1}}{n!}
= \Gamma(q)x^{q-1}(1+x)^{-q};
\end{equation}
express $(1 + x)^{-q}$ as a power series in $x$ at $x = 0$
and use $\Gamma(s) s = \Gamma(s + 1)$.

We use these two formulas (\ref{reviewer1}) and (\ref{reviewer2})
to evaluate
\begin{equation} \label{reviewer3}
I(x) \equiv \frac{1}{2\pi i}\int_{(c)}
  \Gamma(s) \Gamma(s+q-1)\Gamma(1-s) x^{-s} ds.
\end{equation}

Here, we recall that \cite{T2}
\[   
 \frac{1}{2\pi i}\int_{(c)} 
 \mathcal{F}(s) \mathcal{G}(s)x^{-s} ds 
 =
\int_{0}^{\infty} f(z) g(\frac{x}{z})
\frac{dz}{z},
\]
where $\frac{1}{2\pi i}\int_{(c)} \mathcal{F}(s)x^{-s}ds = f(x)$ 
and $\frac{1}{2\pi i}\int_{(c)} \mathcal{G}(s)x^{-s} ds = g(x)$.

We associate $\mathcal{F}$ and $f$ with $\Gamma(s)$ and $e^{-x}$
(by (\ref{reviewer1})), and 
$\mathcal{G}$ and $g$ with $\Gamma(s+q-1)\Gamma(1-s) $
and $\Gamma(q)x^{q-1}(1+x)^{-q}$ (by (\ref{reviewer2})), 
respectively.

This gives
\begin{equation} \label{reviewer4}
I(x) = \Gamma(q)x^{q - 1} 
\int_{0}^{\infty}
e^{-z}
(z+x)^{-q}dz.
\end{equation}

Thus choosing $x \mapsto xpn$, multiplying by 
\[
\Lambda(n)
e^{-2 \pi i p n}p^{-\delta}n^{-\delta - \kappa},
\]
and summing over $n \geq 2$,
we have for $\kappa$ and $\delta$ as in the hypothesis
of the theorem,
\begin{equation} \label{ree3}
\begin{split}
&\frac{1}{2\pi i}\int_{(c)} 
\sum_{n \geq 2} \frac{\Lambda(n) 
e^{- 2 \pi i p n} }
{p^{s + \delta}n^{s + \delta + \kappa}}
\Gamma(s) \Gamma(s+q-1)\Gamma(1-s)x^{-s}ds \\
& = \sum_{n=1}^{\infty} \frac{\Lambda(n)e^{-2\pi i p n}}
{p^{\delta}n^{\delta + \kappa}} I( x pn) \\
& = \Gamma(q) x^{q - 1} \sum_{n=1}^{\infty} \frac{\Lambda(n)
e^{-2 \pi i p n}
}{p^{1 - q + \delta}n^{1-q + \delta + \kappa}} 
\int_{0}^{\infty}   e^{-z}  (z + xpn)^{-q} dz,
\end{split}
\end{equation}
or rearranging,
\[
\begin{split}
& \frac{1}{2\pi i}\int_{(c)} 
\sum_{n \geq 2} \frac{\Lambda(n) e^{- 2\pi i p n}}{p^{s + \delta}
n^{s + \delta + \kappa}}
\Gamma(s) \Gamma(s+q-1)\Gamma(1-s)x^{-s}ds \\
& = \Gamma(q ) x^{-1} \sum_{n=1}^{\infty} 
\frac{\Lambda(n) e^{-2\pi i p n}}{p^{1 - q + \delta}n^{1 - q + \delta + \kappa}} 
\int_{0}^{\infty}   e^{-z}  (z/x + pn)^{-q} dz.
\end{split}
\]

We use the
the change of variables $z = xw$, let $x = 2 \pi e^{\pi i / 2}$, and obtain
\begin{equation} \label{star2}
\begin{split}
& \frac{1}{2\pi i}\int_{(c)} 
\sum_{n \geq 2} \frac{\Lambda(n) e^{- 2\pi i p n} }
{p^{ \delta}n^{s + \delta + \kappa}}
\Gamma(s) \Gamma(s+q-1)\Gamma(1-s)(2 \pi p)^{-s}e^{-\pi i s / 2}ds \\
& = \Gamma(q) \sum_{n=1}^{\infty} 
\frac{\Lambda(n)e^{- 2 \pi i p n}}{p^{1 - q + \delta}n^{1- q + \delta + \kappa}} 
\int_{0}^{\infty}   e^{- 2 \pi i w}  (w + pn)^{-q} dw.
\end{split}
\end{equation}

We choose $q = 1 + \delta$ by analytic continuation; this still
keeps the left side to be in a residue-free form.

Next, we rewrite the integral $\int_{0}^{\infty}$ on the right side
of (\ref{star2}) as
\begin{equation} \label{ppp}
\begin{split}
& \int_{0}^{\infty}  e^{- 2\pi iw} 
 (pn +  w)^{-1 - \delta} dw \\
& =  p^{-\delta } \int_{0}^{\infty}  e^{- 2\pi i pw'} 
 (n +  w')^{-1 - \delta} dw' \\
& = p^{- \delta} \sum_{k \geq 0}\int_{k}^{k + 1}  e^{- 2\pi i pw'} 
 (n +  w')^{-1 - \delta} dw' \\
& =  p^{- \delta} \sum_{k \geq 0}\int_{0}^{1}  e^{- 2 \pi i p (k + r)} 
 (n + k + r)^{-1 - \delta} dr \\
& = p^{- \delta} e^{ 2\pi i p n} \int_{0}^{1} 
\sum_{j \geq n} e^{- 2\pi i p (j + r)} 
 (j + r)^{-1 - \delta} dr,
\end{split}
\end{equation}
where in obtaining the first equality, we made
the change of variables $w = p w'$.

By
\[
\begin{split}
& \sum_{j \geq n} e^{- 2\pi i p (j + r)} 
 (j + r)^{-1 - \delta}  \\
& = \sum_{j \geq 1} e^{- 2\pi i p (j + r)} 
 (j + r)^{-1 - \delta}  - \sum_{1 \leq j < n} e^{- 2\pi i p (j + r)} 
 (j + r)^{-1 - \delta}  \\
& \equiv e^{- 2 \pi i p r}\Phi(p, r, 1 + \delta) 
- \sum_{1 \leq j < n} e^{- 2\pi i p (j + r)} 
 (j + r)^{-1 - \delta} ,
\end{split}
\]
(\ref{ppp}) becomes
\[
\begin{split}
& \int_{0}^{\infty}  e^{- 2\pi iw} 
 (pn +  w)^{-1 - \delta} dw \\
& = p^{- \delta} e^{ 2\pi i p n} \int_{0}^{1} 
\sum_{j \geq n} e^{- 2\pi i p (j + r)} 
 (j + r)^{-1 - \delta} dr \\
& = p^{- \delta} e^{ 2\pi i p n} 
\int_{0}^{1}
\big( e^{- 2 \pi i p r}\Phi(p, r, 1 + \delta)  - \sum_{1 \leq j < n} e^{- 2\pi i p (j + r)} 
 (j + r)^{-1 - \delta} \big) dr,
\end{split}
\]
and so (\ref{star2}) is rewritten as
\begin{equation} \label{star23}
\begin{split}
& \frac{1}{2\pi i}\int_{(c)} 
\sum_{n \geq 2} \frac{\Lambda(n) e^{- 2 \pi i p n} }{p^{ \delta}
n^{s + \delta + \kappa}}
\Gamma(s) \Gamma(s+ \delta)\Gamma(1-s)(2 p \pi)^{-s}e^{-\pi i s / 2}ds \\
& = \Gamma(1 + \delta)
p^{-\delta}
 \sum_{n=2}^{\infty} 
\frac{\Lambda(n)}{n^{\kappa}} \\
& \times 
\int_{0}^{1}
\big( e^{- 2 \pi i p r}\Phi(p, r, 1 + \delta)  - \sum_{1 \leq j < n} e^{- 2 \pi i p (j + r)} 
 (j + r)^{-1 - \delta} \big) dr \\
& \equiv  \Gamma(1 + \delta)
p^{-\delta}(Y_{1} - Y_{2}).
\end{split}
\end{equation}

Here, it is plain that
with the dominated convergence theorem,
\[
Y_{2} = \int_{0}^{1}
 \sum_{n=2}^{\infty} 
\frac{\Lambda(n)}{n^{\kappa}}
 \sum_{1 \leq j < n} e^{-2 \pi i p (j + r)} 
 (j + r)^{-1 - \delta} dr.
\]

In order to analyze $Y_{2}$,
we use a variation of the following relation \cite[pp.60]{T} 
\begin{equation} \label{appc}
\begin{split}
\sum_{j < x} \frac{a_{j}}{j^{s}}
& = \frac{1}{2 \pi i } \int_{c' - iU}^{c' + iU} f(s + w) \frac{x^{w}}{w}dw
+ O(x^{c'} U^{-1} (\sigma + c' - 1)^{-\alpha}) \\
& + O(U^{-1}\psi(2x) x^{1 - \sigma} \log x) + O(U^{-1}\psi(N) x^{1 - \sigma} |x - N|^{-1}), 
\end{split}
\end{equation}
where $x$ is not an integer, $N$ is the integer nearest to $x$,
$c' > 0$, $\sigma + c'  > 1$,
$a_{n} \ll \psi(n)$ for some non-decreasing $\psi$, and the series
\[
f(s) = \sum_{n \geq 1} \frac{a_{n}}{n^{s}}, \quad s = \sigma + it,
\]
converges absolutely for $\sigma > 1$ with
\[
\sum_{n \geq 1} \frac{|a_{n}|}{n^{\sigma}} \ll (\sigma - 1)^{- \alpha}.
\]

In the proof of (\ref{appc}) available in \cite{T},
we replace ``$n$" and ``$x$" by $j + r$ and $n$, respectively,
and obtain (the ratio $(j + r)/n$ is still noninteger for all $n$ and $j$)
\begin{equation} \label{appc4}
\begin{split}
\sum_{j < n} \frac{a_{j}}{(j + r)^{s}}
& = \frac{1}{2 \pi i } \int_{c' - iU}^{c' + iU} f(s + w, r) \frac{n^{w}}{w}dw
+ O(n^{c'} U^{-1} (\sigma + c' - 1)^{-\alpha}) \\
& + O(U^{-1}\psi(2n) n^{1 - \sigma} \log n) + O(U^{-1}\psi(n) n^{1 - \sigma} r^{-1}),
\end{split}
\end{equation} 
where $r \in (0, 1)$ and
\[
f(s, r) \equiv \sum_{n} \frac{a_{n}}{(n + r)^{s}}.
\]

Choosing $s \mapsto 1 + \delta$ and
$a_{j} \mapsto e^{- 2 \pi i p j}$ in (\ref{appc4}),
we have
\[
\begin{split}
& \sum_{j < n} \frac{e^{- 2\pi i p j}}{(j + r)^{1 + \delta}} \\
& = \frac{1}{2 \pi i } \int_{c' - iU}^{c' + iU} 
\Phi(p, r, 1 + \delta + w)
\frac{n^{w}}{w}dw  + O(U^{-1} (\log n)^{2}  (n + r)^{c'}r^{-1}  (\sigma + c' - 1)^{-\alpha}).
\end{split}
\]
Furthermore, multiplying both sides by $\Lambda(n)n^{- \kappa}$
and summing all over the positive integers $n \geq 2$, we get
\begin{equation} \label{45dd}
\begin{split}
& \sum_{n \geq 2} \frac{\Lambda(n)}{n^{\kappa}} 
\sum_{1 \leq j < n} \frac{e^{- 2\pi i p j}}{(j + r)^{1 + \delta}} \\
& = \frac{1}{2 \pi i } \int_{c' - iU}^{c' + iU}- \frac{\zeta'}{\zeta}
(\kappa - w) \frac{\Phi(p, r, 1 + \delta + w) }{w}dw \\
& + O_{r}(U^{-1} H (\sigma + c' - 1)^{-\alpha}),
\end{split}
\end{equation}
where
\[
H \equiv \sum_{n \geq 2}
(\log n)^{2} \Lambda(n)  n^{-\kappa + c'}.
\]
By (\ref{45dd}), we see that
\begin{equation} \label{shikaku}
\begin{split}
 \frac{1}{2 \pi i } 
& \int_{c' - iU}^{c' + iU}- \frac{\zeta'}{\zeta}
(\kappa - w) \frac{ \Phi(p, r, 1 + \delta + w) }{w}dw \\
& \to 
 \sum_{n \geq 2} \frac{\Lambda(n)}{n^{\kappa}} 
\sum_{1 \leq j < n} \frac{e^{- 2\pi i p j}}{(j + r)^{1 + \delta}}
\end{split}
\end{equation}
uniformly in $r \in [1/N, 1 - 1/N]$
for any fixed $N$ as $U \to \infty$.

Hence, the integral in $Y_{2}$
is rewritten as
\begin{equation} \label{rewwef}
\begin{split}
& \lim_{N \to \infty}\int_{1/N}^{1 - 1/N}
\sum_{n \geq 2} \frac{\Lambda(n)}{n^{\kappa}} 
\sum_{1 \leq j < n} \frac{e^{- 2 \pi i p j}}{(j + r)^{1 + \delta}}
e^{- 2\pi i p r}
  dr \\
& = \lim_{N  \to \infty} \int_{1/N}^{1 - 1/N} 
 \frac{1}{2 \pi i } \int_{c' - i \infty}^{c' + i\infty}- \frac{\zeta'}{\zeta}
(\kappa - w) \frac{\Phi(p, r, 1 + \delta + w) }{w}dw
e^{- 2\pi i p r} dr \\
& =  \lim_{N  \to \infty} 
 \frac{1}{2 \pi i } \int_{c' - i \infty}^{c' + i\infty}
- \frac{\zeta'}{\zeta}
(\kappa - w) 
\int_{1/N}^{1 - 1/N} \Phi(p, r, 1 + \delta + w)
e^{- 2\pi i p r} dr 
\frac{dw}{w};
\end{split}
\end{equation}
the first equality is by pointwise convergence of (\ref{shikaku}), and
the second by uniform convergence of (\ref{shikaku}) for each $N$.

But if $c'$ and $\kappa$ are sufficiently large so that $\kappa - c' > 1$ 
and $1 + \delta + c' > 2$ (as in the hypothesis of the theorem), then
with integration by parts, it is easy to show that
\begin{equation} \label{bdh}
\int_{1/N}^{1 - 1/N} \Phi(p, r, 1 + \delta + w) e^{-2 \pi i p r} dr \ll (|t| + 1)^{-1}, \quad w = c' + it,
\end{equation}
and so the last integral $\int_{(c')}$ in (\ref{rewwef}) converges absolutely.
This in turn enables us to put the limit $N \to \infty$
inside the integral symbol $\int_{(c')}$ (use the 
dominated convergence theorem); we get 
\begin{equation} \label{rewwefe}
Y_{2} =   
 \frac{1}{2 \pi i } \int_{c' - i \infty}^{c' + i\infty}
- \frac{\zeta'}{\zeta}
(\kappa - w) 
\int_{0}^{1}
\Phi(p, r, 1 + \delta + w)
 e^{- 2\pi i p r} dr 
\frac{dw}{w}.
\end{equation}

By (\ref{rewwefe}),
(\ref{star23}) becomes 
\begin{equation} \label{star24}
\begin{split}
& \frac{1}{2\pi i}\int_{(c)} 
\sum_{n \geq 2} \frac{\Lambda(n) e^{-2 \pi i p n} }{
n^{s + \delta + \kappa}}
\Gamma(s) \Gamma(s+ \delta)\Gamma(1-s) (2 p \pi)^{-s}e^{-\pi i s / 2}ds \\
& = \Gamma(1 + \delta)
 \sum_{n=1}^{\infty} 
\frac{\Lambda(n)}{n^{\kappa}} \int_{0}^{1}
 e^{- 2 \pi i p r}\Phi(p, r, 1 + \delta)dr \\
& -  \frac{\Gamma(1 + \delta)}{2 \pi i }
\int_{c' - i \infty}^{c' + i\infty}
- \frac{\zeta'}{\zeta}
(\kappa - w)
\int_{0}^{1}
\Phi(p, r, 1 + \delta + w)
 e^{- 2\pi i p r} dr 
\frac{dw}{w}.
\end{split}
\end{equation}

Finally, we rewrite for $0 < \text{Re}(s) < 1$,
\[
\begin{split}
\int_{0}^{1}
\Phi(p, r, s)
 e^{- 2\pi i p r} dr  
& = \sum_{n \geq 1} \int_{0}^{1}
 e^{- 2\pi i p(n + r)} 
(n + r)^{-s}dr \\
& =  \sum_{n \geq 1} \int_{n}^{n + 1}
 e^{- 2\pi i p y} 
y^{-s}dy \\
& =  \int_{1}^{\infty}
 e^{- 2\pi i p y} 
y^{-s}dy \\
& = \left( \int_{0}^{\infty} - \int_{0}^{1} \right)
 e^{- 2\pi i p y} 
y^{-s}dy \\
& \equiv J(p, s)+ K(p, s).
\end{split}
\]

The term $J(p, s)$ is rewritten with the residue theorem and change of variables as
\[
\begin{split}
J(p, s) = \int_{0}^{\infty} 
 e^{- 2\pi i p y} 
y^{-s}dy
& = \int_{0}^{\infty} 
 e^{- 2\pi p z} 
[(-i)z]^{-s} (-i)dz \\
& = [(-i)(2\pi p)]^{-s + 1} \int_{0}^{\infty}
 e^{- w} 
w^{1 - s - 1 } dz \\
& =  [(-i)(2\pi p)]^{-s + 1} \Gamma(1 - s) \\
& = (2 \pi p)^{s - 1}  e^{(s -1 ) \pi i / 2}\Gamma(1 - s),
\end{split}
\]
which is readily extended to the region outside the strip
$0 < \text{Re}(s) < 1$.

With this, (\ref{star24}) becomes,
\begin{equation} \label{star25}
\begin{split}
& \frac{1}{2\pi i}\int_{(c)} 
\sum_{n \geq 2} \frac{\Lambda(n) e^{-2 \pi i p n} }{
n^{s + \delta + \kappa}}
\Gamma(s) \Gamma(s+ \delta)\Gamma(1-s) (2 p \pi)^{-s}e^{-\pi i s / 2}ds \\
& = \Gamma(1 + \delta)
(- \frac{\zeta'}{\zeta}(\kappa))
[ (2 \pi p)^{-\delta}  e^{\delta \pi i / 2}\Gamma( - \delta)
+ K(p, 1 + \delta)]
\\
& -  \frac{\Gamma(1 + \delta)}{2 \pi i }
\int_{c' - i \infty}^{c' + i\infty}
- \frac{\zeta'}{\zeta}
(\kappa - w) \\
& \times 
 [(2 \pi p)^{\delta + w}  e^{(\delta + w ) \pi i / 2}\Gamma( - \delta - w)
+ K(p, 1 + \delta + w)]
\frac{dw}{w}.
\end{split}
\end{equation}

This completes the proof of Theorem \ref{main1}.

\subsection{The second proof}
The second proof
avoids algebraic manipulations appearing in the first proof, and so deserves to be 
shown here.\footnote{The main part of the argument was 
brought by Gemini.}

Let all the variables have the same meanings as stated in Theorem \ref{main1}.

We prove the theorem in the form (\ref{star24}).

By the definition of $\Phi(p, r, 1 + \delta)$, we rewrite the first integral
on the right side of (\ref{star24}) as
\begin{equation} \label{star241}
\begin{split}
 \int_{0}^{1}
 e^{- 2 \pi i p r}\Phi(p, r, 1 + \delta)dr 
& =  \int_{0}^{1} e^{-2 \pi i p r} \sum_{j \geq 1} e^{- 2 \pi i p j} (j + r)^{-1 - \delta} dr \\
& = \int_{0}^{1} \sum_{j \geq 1} e^{- 2 \pi i p (j + r)} (j + r)^{-1 - \delta} dr \\
& = \int_{1}^{\infty} e^{- 2 \pi i p x} x^{-1 - \delta} dx.
\end{split}
\end{equation}

Since $\kappa - c' > 1$, we can rewrite the factor $- (\zeta'/\zeta)(\kappa - w)$
inside the second integral of the right side of (\ref{star24}) as
\[
- \frac{\zeta'}{\zeta}(\kappa - w)
= \sum_{n \geq 1} \frac{\Lambda(n)}{n^{\kappa - w}}.
\]
Hence, with (\ref{star241}), we can rewrite the second integral of the right side of (\ref{star24}) as
\begin{equation} \label{star242}
\begin{split}
& \frac{1}{2 \pi i} \int_{c' - i \infty}^{c' + i\infty}
- \frac{\zeta'}{\zeta}
(\kappa - w)
\int_{0}^{1}
\Phi(p, r, 1 + \delta + w)
 e^{- 2\pi i p r} dr 
\frac{dw}{w} \\
& =  \frac{1}{2 \pi i}\int_{c' - i \infty}^{c' + i\infty}
 \sum_{n \geq 1} \frac{\Lambda(n)}{n^{\kappa - w}}
 \int_{1}^{\infty} e^{- 2 \pi i p x} x^{-1 - \delta - w} dx
\frac{dw}{w}
\\
& =  \frac{1}{2 \pi i}\int_{c' - i \infty}^{c' + i\infty}
 \sum_{n \geq 1} \frac{\Lambda(n)}{n^{\kappa}}
 \int_{1}^{\infty} e^{- 2 \pi i p x} x^{-1 - \delta} (n/x)^{ w} dx
\frac{dw}{w}.
\end{split}
\end{equation}
By $\kappa - c' > 1$ and using integration by parts to the inner integral,
it is plain that we can interchange the summation over $n$ and the integration in the variable $w$. 
We get
\begin{equation} \label{star243}
\begin{split}
& \frac{1}{2 \pi i} \int_{c' - i \infty}^{c' + i\infty}
- \frac{\zeta'}{\zeta}
(\kappa - w)
\int_{0}^{1}
\Phi(p, r, 1 + \delta + w)
 e^{- 2\pi i p r} dr 
\frac{dw}{w} \\
& = \frac{1}{2 \pi i} \sum_{n \geq 1} \frac{\Lambda(n)}{n^{\kappa}}
 \int_{c' - i \infty}^{c' + i\infty}
 \int_{1}^{\infty} e^{- 2 \pi i p x} x^{-1 - \delta} (n/x)^{ w} dx
\frac{dw}{w}.
\end{split}
\end{equation}

It is well-known \cite{T} that 
\begin{equation} \label{PR}
 \frac{1}{2 \pi i} \int_{c' - i T}^{c' + iT} (n/x)^{ w} 
\frac{dw}{w}
=
\begin{cases}
1 + O\left(\frac{(n/x)^{c'}}{T \log (n/x)} \right) &  x < n  \\
O\left(\frac{(n/x)^{c'}}{T \log (n/x)} \right) &  x > n.
\end{cases}
\end{equation}
Hence, for each fixed $n$, the convergence as $T \to \infty$ is uniform, and so we can
interchange the integrals of (\ref{star243}). We obtain
\[
\begin{split}
& \frac{1}{2 \pi i} \int_{c' - i \infty}^{c' + i\infty}
- \frac{\zeta'}{\zeta}
(\kappa - w)
\int_{0}^{1}
\Phi(p, r, 1 + \delta + w)
 e^{- 2\pi i p r} dr 
\frac{dw}{w} \\
& = \frac{1}{2 \pi i} \sum_{n \geq 1} \frac{\Lambda(n)}{n^{\kappa}}
 \int_{1}^{\infty}e^{- 2 \pi i p x} x^{-1 - \delta} \int_{c' - i \infty}^{c' + i\infty}
  (n/x)^{ w}
\frac{dw}{w} dx,
\end{split}
\]
or with (\ref{PR}), 
\begin{equation} \label{star244}
\begin{split}
& \frac{1}{2 \pi i} \int_{c' - i \infty}^{c' + i\infty}
- \frac{\zeta'}{\zeta}
(\kappa - w)
\int_{0}^{1}
\Phi(p, r, 1 + \delta + w)
 e^{- 2\pi i p r} dr 
\frac{dw}{w} \\
& = \sum_{n \geq 1} \frac{\Lambda(n)}{n^{\kappa}}
 \int_{1}^{n}e^{- 2 \pi i p x} x^{-1 - \delta} dx.
\end{split}
\end{equation}

Recalling (\ref{star241}), the right side of (\ref{star24}) is, with (\ref{star244}),
\begin{equation} \label{rhs24}
(\text{the right side of (\ref{star24})})
= \Gamma(1 + \delta) \sum_{n \geq 1} \frac{\Lambda(n)}{n^{\kappa}}
 \int_{n}^{\infty}e^{- 2 \pi i p x} x^{-1 - \delta} dx.
\end{equation}

Next, we transform the $n$th term of (\ref{rhs24}) further to match that of (\ref{star24}), as follows.

With the change of variables $x = u + n$, we get
\[
\begin{split}
 \frac{\Lambda(n)}{n^{\kappa}}
 \int_{n}^{\infty}e^{- 2 \pi i p x} x^{-1 - \delta} dx
=  \frac{\Lambda(n) e^{- 2 \pi i pn}}{n^{\kappa}}
 \int_{0}^{\infty}e^{- 2 \pi i p u} (u + n)^{-1 - \delta} du,
\end{split}
\]
and with $u = nv$,
\begin{equation} \label{rhs241}
\begin{split}
 \frac{\Lambda(n)}{n^{\kappa}}
 \int_{n}^{\infty}e^{- 2 \pi i p x} x^{-1 - \delta} dx
=  \frac{\Lambda(n) e^{- 2 \pi i pn}}{n^{\kappa + \delta}}
 \int_{0}^{\infty}e^{- 2 \pi i p nv} (v + 1)^{-1 - \delta} dv.
\end{split}
\end{equation}

On the other hand, the Mellin transform of 
\[
f(\alpha) \equiv  \int_{0}^{\infty}e^{- i \alpha v} (v + 1)^{-1 - \delta} dv
\]
gives 
\begin{equation} \label{lhs21}
\begin{split}
 M[f](s) & = 
\int_{0}^{\infty}  
\left[ \int_{0}^{\infty}e^{- i \alpha v} (v + 1)^{-1 - \delta} dv \right]
\alpha^{s - 1} d\alpha \\
& = \int_{0}^{\infty}
\left[ \int_{0}^{\infty}  
e^{- i \alpha v}
\alpha^{s - 1} d\alpha
\right]
 (v + 1)^{-1 - \delta} dv.
\end{split}
\end{equation}
The inner integral is rewritten as
\[
\int_{0}^{\infty}  
e^{- i \alpha v}
\alpha^{s} \frac{d\alpha}{ \alpha}
=  \int_{0}^{\infty}  
e^{- i u}
(u/v)^{s} \frac{du}{u}
= v^{-s} \int_{0}^{\infty}  
e^{- i u} u^{s - 1}du,
\]
and so (\ref{lhs21}) becomes
\begin{equation} \label{lhs22}
\begin{split}
 M[f](s)
& = \int_{0}^{\infty}
\left[ \int_{0}^{\infty}  
e^{- i \alpha v}
\alpha^{s - 1} d\alpha
\right]
 (v + 1)^{-1 - \delta} dv \\
& = \int_{0}^{\infty}
e^{- i u}
u^{s - 1} du
\int_{0}^{\infty}  v^{-s} (v + 1)^{-1 - \delta} dv.
\end{split}
\end{equation}
Now, in the formula
\[
\int_{0}^{\infty} e^{- at}t^{s - 1} dt = a^{-s} \Gamma(s),
\]
we let $a \to i = e^{\pi i / 2}$ and obtain
\begin{equation} \label{Gammasan}
\int_{0}^{\infty} e^{- it}t^{s - 1} dt = e^{- \pi i s /2} \Gamma(s).
\end{equation}

We recall \cite{WW} that the beta function $B(a, b)$ satisfies
\[
B(a, b) = \int_{0}^{\infty} \frac{x^{a -  1}}{(1 + x)^{a + b}} dx
= \frac{\Gamma(a) \Gamma(b)}{\Gamma(a + b)}.
\]
Choosing $a = 1 - s$ and $b = s + \delta$, we find out that
\begin{equation} \label{betaa}
\int_{0}^{\infty} \frac{x^{-s}}{(1 + x)^{1 + \delta}} dx
= \frac{\Gamma(1 - s) \Gamma(s + \delta)}{\Gamma(1 + \delta)}.
\end{equation}

Inserting (\ref{Gammasan}) and (\ref{betaa}) into (\ref{lhs22}), we get
\begin{equation} \label{lhs23}
\begin{split}
 M[f](s)
&  = \int_{0}^{\infty}
e^{- i u}
u^{s - 1} du
\int_{0}^{\infty}  v^{-s} (v + 1)^{-1 - \delta} dv \\
& = 
 \frac{ e^{- \pi i s /2} \Gamma(s)\Gamma(1 - s) \Gamma(s + \delta)}{\Gamma(1 + \delta)}.
\end{split}
\end{equation}
Applying the Mellin inversion formula, we find out that
\begin{equation} \label{lhs24}
\begin{split}
\Gamma(1 + \delta)f(\alpha)
  = \frac{1}{2 \pi i}
\int_{(c')}  e^{- \pi i s /2} \Gamma(s)\Gamma(1 - s) \Gamma(s + \delta)
\alpha^{-s} ds.
\end{split}
\end{equation}
Substituting $\alpha = 2 \pi p n$ in (\ref{lhs24}) and
then multiplying by $\Lambda(n)n^{- \delta - \kappa}$, we see that
\begin{equation} \label{lhs25}
\begin{split}
&  
\frac{\Gamma(1 + \delta)\Lambda(n)}{n^{\delta + \kappa}}
 \int_{0}^{\infty}e^{- 2 \pi i  p  n v} (v + 1)^{-1 - \delta} dv \\
&  = \frac{1}{2 \pi i}
\int_{(c')}
\frac{\Lambda(n)}{n^{\delta + \kappa + s}}
  e^{- \pi i s /2} \Gamma(s)\Gamma(1 - s) \Gamma(s + \delta)
(2 \pi p)^{-s} ds,
\end{split}
\end{equation}
which shows that, with (\ref{rhs241}), the $n$th term of (\ref{rhs24})
matches that of the left side of (\ref{star24}).

This completes the proof of Theorem \ref{main1}.

\section{
Proof of
Theorem \ref{supraclaim}
} 

In this section, we give a proof of Theorem \ref{supraclaim}.

We use the same notations appearing in the statement of the theorem. 

First, we transform the left side of (\ref{thm1}) into that of Theorem \ref{supraclaim}.

As we make $\kappa = 0$ and $\nu > 0$ small, we
shift the path of the integral to the vertical line $\text{Re}(s) = -1/2$. 
Recalling Lemma \ref{rwri}, we obtain for $p = a/b$, 
\begin{equation} \label{w1}
\begin{split}
& \frac{1}{2\pi i}\int_{(c)} 
\sum_{n \geq 2} \frac{\Lambda(n) e^{-2 \pi i p n} }{
n^{s + \delta}}
\Gamma(s) \Gamma(s+ \delta)\Gamma(1-s) (2 p\pi)^{-s}e^{-\pi i s / 2}ds \\
& =- \sum_{\chi }\sum_{\rho_{\chi}}A(a, b; \chi) (2 p \pi)^{- \rho_{\chi} + \delta}
e^{-\pi i (\rho_{\chi} - \delta ) / 2}  \Gamma(
\rho_{\chi} - \delta
)
\Gamma(\rho_{\chi})\Gamma(1- \rho_{\chi} + \delta) \\
& + M(0, a/b) \Gamma(\delta) 
- \sum_{\chi}A(a, b; \chi) \text{res}_{s = - \delta}\left[ \frac{L'}{L}(s + \delta) 
\Gamma(s) \Gamma(s+ \delta) \Gamma(1-s) (2 p\pi)^{-s}e^{-\pi i s / 2}\right] \\
& - \sum_{\chi^{*}} \sum_{n}A(a, b; \chi^{*}) 
 \Gamma( h^{*}(n) - \delta) \Gamma(h^{*}(n)) 
\Gamma(1-h^{*}(n) + \delta) (2 p\pi)^{- h^{*}(n) + \delta}e^{-\pi i ( h^{*}(n) - \delta) / 2} \\
& + \text{res}_{s = - \delta}[E_{M}(s + \delta, a/b)\Gamma(s) \Gamma(s+ \delta) \Gamma(1-s) (2 p\pi)^{-s}e^{-\pi i s / 2} ] \\
&+ \frac{1}{2\pi i}\int_{(-1/2)} 
M(s + \delta, a/b)
\Gamma(s) \Gamma(s+ \delta)\Gamma(1-s) (2 p\pi)^{-s}e^{-\pi i s / 2}ds.
\end{split}
\end{equation}

Second, using residue theorem, we attempt to let $\kappa = 0$ and $\nu > 0$ be small in the right side of (\ref{thm1})
as follows.

On the right side of (\ref{thm1}), we apply the residue theorem 
to the second integral and get
\begin{equation} \label{thm1rhs11}
\begin{split}
& \frac{1}{2 \pi i } \int_{c' - i \infty}^{c' + i\infty}
- \frac{\zeta'}{\zeta}
(\kappa - w) 
[
J(p, 1 + \delta + w) + K(p, 1 + \delta + w)
]
\frac{dw}{w} \\
& =
- \sum_{\rho_{j}}
(\kappa - \rho_{j})^{-1}
[J(p, 1 + \delta + \kappa - \rho_{j}) + K(p,  1 + \delta + \kappa - \rho_{j})]
\\
&
 + (\kappa - 1)^{-1} 
[
J(p, \delta + \kappa) + K(p, \delta + \kappa)
]   
\\
& + \frac{1}{2 \pi i } \int_{2 - h - i \infty}^{2 - h + i\infty}
- \frac{\zeta'}{\zeta}
(\kappa - w) 
[J(p, 1 + \delta + w) + K(p, 1 + \delta + w)]
\frac{dw}{w},
\end{split}
\end{equation}
with $h > 0$ arbitrary.

With (\ref{gf}), it is easy to see that
\[
 (\kappa - \rho_{j})^{-1}(2\pi p)^{\delta + \kappa - \rho_{j}} 
e^{(\delta + \kappa - \rho_{j})\pi i / 2} \Gamma(-\delta - \kappa + \rho_{j})  \ll \gamma_j^{\eta_{j} - \nu - 1}.
\]
For $\sigma = 1/2 + \eta$ ($0 \leq \eta < 1/2$), we have \cite{T}
\[
N(1/2 + \eta, T) \ll  T^{1 - \eta} (\log T)^5.
\]
Writing $\gamma_j^{\eta_{j}} = 1 + (\log \gamma_{j}) \int_{0}^{\eta_j} \gamma_j^{\eta}\, d\eta$ and expressing the sum over $\gamma_j > 0$ via Riemann--Stieltjes integration with respect to $N(1/2 + \eta, t)$, we obtain by integration by parts in $t$
\begin{equation} \label{otto}
\begin{split}
 \sum_{\text{Re}(\rho_{j}) \geq 1/2} \gamma_j^{\eta_{j} - \delta - 1} 
& \ll \int_{1}^{\infty} t^{-\nu - 1}\, dN(1/2, t) + \int_{0}^{1/2} \left( \int_{1}^{\infty} 
t^{\eta - \nu - 1} (\log t)\, d_t N(1/2 + \eta, t) \right) d\eta  \\
&\ll \int_{1}^{\infty} t^{-\nu - 2} N(1/2, t)\, dt + \int_{0}^{1/2} \int_{1}^{\infty} t^{\eta - \nu - 2} (\log t) N(1/2 + \eta, t)\, dt\, d\eta  \\
&\ll \int_{1}^{\infty} t^{-\nu - 1} (\log t)^6 \, dt < \infty.
\end{split}
\end{equation}
Thus, the first series on the right-hand side of (\ref{thm1rhs11}) converges absolutely and uniformly on compact subsets of 
$\text{Re}(\delta) > 0$.

Combining (\ref{thm1}), (\ref{w1}), (\ref{thm1rhs11}), and (\ref{otto}) together,
we obtain the first claim of Theorem \ref{supraclaim}.

Next, the continuity property of the right members of (\ref{thm1})
is proved as follows. By applying Weierstrass $M$-test to
the first series on the right side of (\ref{thm1rhs11})  (cf. (\ref{otto})),
it is easy to see that the series converges uniformly in $p$. 
Then, we know that {\it if $f_{n}(x)$ is continuous and uniformly
convergent to $f(x)$, then $f(x)$ is continuous}.
Hence,
the first series on the right side of (\ref{thm1rhs11}) 
 is continuous in $p \in (0, 1)$ for each $\nu > 0$.

This completes the proof of
Theorem \ref{supraclaim}.

\section{
Theorem \ref{supraclaim}
over a short interval
}

In this section, we consider the case of Theorem \ref{supraclaim}
with large $b$.

For convenience, we choose $b = \mathfrak{b}$ to be prime. Then all the characters 
appearing in the theorem become primitives, and the error term $E_{M}(s, a/b)$
of Lemma \ref{rwri} is simplified to
\[
E_{M}(s, a/ \mathfrak{b})
= (\log \mathfrak{b})\frac{ \mathfrak{b}^{-s}}{1- \mathfrak{b}^{-s}}.
\]

Besides, we can give a bound as $\mathfrak{b} \to \infty$ for the integral $E_{4}$.
\begin{lem} \label{boundE4}
Let $\delta = \nu + iT$ and $\nu > 0$ small.
Let $\mathfrak{b}$ be prime. Then as $\mathfrak{b} \to \infty$, we have
\[
\begin{split}
E_{4} & =  \frac{1}{2\pi i}\int_{(-1/2)} 
\left[
\sum_{\chi}A(a, \mathfrak{b}; \chi) O(\log[(s + \delta)(1 - s - \delta)])
+ O(\log \mathfrak{b})
\right] \\
& \times \Gamma(s) \Gamma(s+ \delta)\Gamma(1-s) (2 p\pi)^{-s}e^{-\pi i s / 2}ds.
\end{split}
\]
\end{lem}

\begin{proof}

By taking the logarithmic derivative of the functional equations for $L(s, \chi)$ \cite{H}
\begin{equation} \label{feL1}
\begin{split}
& \pi^{-(1 - s)/2} q^{(1 - s)/2} \Gamma[(1 - s)/2] L(1 - s, \overline{\chi}) \\
& = \frac{q^{1/2}}{\tau(\chi)} \pi^{-s/2} q^{s/2} \Gamma(s/2) L(s, \chi), \quad
\chi(-1) = 1
\end{split}
\end{equation}
and
\begin{equation} \label{feL2}
\begin{split}
& \pi^{-(2 - s)/2} q^{(2 - s)/2} \Gamma[(2 - s)/2] L(1 - s, \overline{\chi}) \\
& = \frac{iq^{1/2}}{\tau(\chi)} \pi^{-(1 + s)/2} q^{(1 + s)/2} \Gamma[(1 +s)/2)] L(s, \chi), \quad
\chi(-1) = -1,
\end{split}
\end{equation}
it is easy to show that 
\begin{equation} \label{feL1.2}
\begin{split}
\log \frac{\pi}{q} - \frac{1}{2} \frac{\Gamma'}{\Gamma}[(1 - s)/2] 
- \frac{1}{2}\frac{\Gamma'}{\Gamma}(s/2)
- \frac{L'}{L}(1 - s, \overline{\chi})
=   \frac{L'}{L}(s, \chi),
\quad \chi(-1) = 1
\end{split}
\end{equation}
and
\begin{equation} \label{feL2.2}
\begin{split}
\log \frac{\pi}{q} - \frac{1}{2} \frac{\Gamma'}{\Gamma}[(2 - s)/2] 
- \frac{1}{2}\frac{\Gamma'}{\Gamma}[(1 + s)/2]
- \frac{L'}{L}(1 - s, \overline{\chi})
=   \frac{L'}{L}(s, \chi), \quad
\chi(-1) = -1.
\end{split}
\end{equation}
With \cite{H}
\[
\frac{\Gamma'}{\Gamma}(s) = \log s + O(|s|^{-1}), \quad -\pi + \epsilon < \arg s < \pi - \epsilon,  \quad
\epsilon > 0,
\quad |s| \to \infty,
\]
(\ref{feL1.2}) and (\ref{feL2.2}) can be put together in the simpler form
\begin{equation} \label{feLLL}
\log \frac{\pi}{q} + O(\log(s(1 - s)))
 - \frac{L'}{L}(1- s, \overline{\chi}) = \frac{L'}{L}(s, \chi), \quad \forall \chi.
\end{equation}

Substituting (\ref{feLLL}) into Lemma \ref{rwri}, $E_{4}$ is estimated as
\[
\begin{split}
 E_{4} & = \frac{1}{2\pi i}\int_{(-1/2)} 
M(s + \delta, a/\mathfrak{b})
\Gamma(s) \Gamma(s+ \delta)\Gamma(1-s) (2 p\pi)^{-s}e^{-\pi i s / 2}ds \\
& = \frac{1}{2\pi i}\int_{(-1/2)} 
\left[
- \sum_{\chi}A(a, \mathfrak{b}; \chi) \frac{L'}{L}(s + \delta, \chi)
+ E_{M}(s + \delta, a/\mathfrak{b})
\right] \\
& \times \Gamma(s) \Gamma(s+ \delta)\Gamma(1-s) (2 p\pi)^{-s}e^{-\pi i s / 2}ds \\
& =  \frac{1}{2\pi i}\int_{(-1/2)} 
\left[
\sum_{\chi}A(a, \mathfrak{b}; \chi)
\left( \frac{L'}{L}(1 - s - \delta, \overline{\chi})
+ O(   \log[(s + \delta)(1 - s - \delta)])
\right)
+ O(\log \mathfrak{b})
\right] \\
& \times
\Gamma(s) \Gamma(s+ \delta)\Gamma(1-s) (2 p\pi)^{-s}e^{-\pi i s / 2}ds.
\end{split} 
\]

As proceeded in (\ref{MtoL}) above,
we can recognize the inner sum as, by $1 - \text{Re}(s) - \nu > 1$,
\[
\begin{split}
\sum_{\chi}A(a, \mathfrak{b}; \chi) \frac{L'}{L}(1 - s - \delta, \overline{\chi})
& = - \sum_{\chi}  A(a, \mathfrak{b}; \chi)  \sum_{n \geq 1} 
\frac{
\Lambda(n)  \overline{\chi}(n)
}{n^{1 - s - \delta}} \\
& = - \sum_{n \geq 1} \frac{\Lambda(n)}{n^{s}}
 \sum_{1 \leq x \leq \mathfrak{b}} e^{-2 \pi i a x/\mathfrak{b}} \sum_{\chi} \frac{\overline{\chi}(xn)}{\phi(\mathfrak{b})}.
\end{split}
\]

By the well-known identity
\[
\sum_{\chi} \chi(xn) = 
\begin{cases}
1 & xn \equiv 1  \pmod{\mathfrak{b}} \\
0 & \text{otherwise},
\end{cases}
\]
if we denote $\hat{n}n \equiv 1 \pmod{\mathfrak{b}}$, then we find out that
\[
\begin{split}
\sum_{\chi}A(a, \mathfrak{b}; \chi) \frac{L'}{L}(1 - s - \delta, \overline{\chi})
& = - \sum_{n \geq 1} \frac{\Lambda(n)}{n^{1 - s - \delta}}
 \sum_{1 \leq x \leq \mathfrak{b}} e^{-2 \pi i a x/\mathfrak{b}} \sum_{\chi} \frac{\overline{\chi}(xn)}{\phi(\mathfrak{b})} \\
& =  - \sum_{n \geq 1} \frac{\Lambda(n)}{n^{1 - s - \delta}} e^{-2 \pi i a \hat{n}/\mathfrak{b}}
+ \sum_{l \geq 1} (\log \mathfrak{b}) \mathfrak{b}^{-(1 - s - \delta)l} e^{-2\pi i a \hat{ \mathfrak{b}}^{l} / \mathfrak{b}}.
\end{split}
\]
This expression is bounded for any sequence of rationals $a/\mathfrak{b}$ which approaches any real number in $(0, 1)$.

This completes the proof of the lemma.

\end{proof}

At this point, let us pose a question. By Corollary \ref{claim}, there exists
a rational number $p_{0} \in (0, 1)$ for which the total of the left members
of Theorem \ref{supraclaim} is large, being associated with $\eta_{n_{1}, p} = \text{Re}(\rho_{n_{1}, p}) - 1/2$.

What could we say in analyzing the values of the left members of the theorem for $p$ in a short interval
$(p_{0}, p_{0} + \epsilon)$? As $\epsilon > 0$ becomes small, the denominator $\mathfrak{b}$ of
the number $p = a/\mathfrak{b} \in (p_{0}, p_{0} + \epsilon)$ would naturally become large.
In general, there is no guarantee that as $\mathfrak{b}$ increases,
all the terms of the left members of the theorem would be uniformly bounded;
as $\mathfrak{b} \to \infty$, the number of $L(s, \chi)$'s to consider would increase. 

Nevertheless, with the continuity claim of the theorem, 
the following is obvious.
\begin{thm} \label{princ1}
As $a/\mathfrak{b} \to p_{0}$,
the total of the left members of Theorem \ref{supraclaim}
converges to that evaluated at $p_{0}$.
\end{thm}

From our analysis presented heretofore, we conclude this paper by the following
general principle: \em{the nontrivial zeros of Dirichlet $L$-functions
are distributed in order for the left members of Theorem \ref{supraclaim}
to be continuous in $p$.}


\begin{thebibliography}{0}
\bibitem{H} 
H. Davenport, 
Multiplicative Number Theory,
3rd ed.,
Springer, 
2000.

\bibitem{HL}
G. H. Hardy and J. E. Littlewood, {\em
Some problems of `Partitio numerorum'; \Romannum{3}: On
the expression of a number as a sum of
primes,
} Acta. Math. {\bf 44} (1923)


\bibitem{008}
M. Jutila, 
Lectures on a Method in the Theory of Exponential Sums,
Tata Inst. Fund. Res. Lectures on Math. and Phys.., 80,
Published for the Tata Institute of Fundamental 
Research, Bombay, by Springer-Verlag, Berlin, 1987.

\bibitem{L} 
S. Lang,
Complex Analysis,
3rd ed., 
Springer,
1993.

\bibitem{T2} 
E. C. Titchmarsh,
Introduction to the Theory of Fourier Integrals,
2nd ed., 
Oxford, 
1948.

\bibitem{T} 
E. C. Titchmarsh,
The Theory of the Riemann Zeta-function, 
2nd ed., 
Oxford, 1986.

\bibitem{WW}
E. T. Whittaker and G. N. Watson,
A Course of Modern Analysis, 3rd ed., Cambridge
University Press, Cambridge, 1920.

\end{thebibliography}
\end{document}